\documentclass{article}
\usepackage[preprint]{neurips_2025}
\usepackage[utf8]{inputenc} 
\usepackage[T1]{fontenc}    
\usepackage{hyperref}       
\usepackage{url}            
\usepackage{booktabs}       
\usepackage{amsfonts}       
\usepackage{nicefrac}       
\usepackage{microtype}      
\usepackage{xcolor}         
\usepackage{amsmath,amsfonts}
\usepackage{algorithm}
\usepackage{algorithmic}

\title{Using a generative model for out-of-sample testing of two-stage stochastic programs}
\workshoptitle{ML×OR: Mathematical Foundations and Operational Integration of Machine Learning for Uncertainty-Aware Decision-Making}

\author{Ashutosh Shukla, John J. Hasenbein, Erhan Kutanoglu \\
  Graduate Program in Operations Research and Industrial Engineering\\
  The University of Texas at Austin\\
  Austin, Texas 78712 \\
  \texttt{\{ashutosh.shukla@utexas.edu, has@me.utexas.edu, erhank@austin.utexas.edu\}} \\
}

\begin{document}
\maketitle
\begin{abstract}
Stochastic programming models for decision-making under uncertainty often suffer from scenario scarcity, where obtaining representative samples of uncertain parameters requires expensive simulations or measurements. This work presents a framework that leverages the Normal-to-Anything (NORTA) generative model to enhance the reliability of two-stage stochastic programming solutions through comprehensive out-of-sample testing when scenario data is limited. The NORTA model efficiently generates synthetic scenarios that preserve both marginal distributions and correlation structures from limited available data, offering a computationally tractable alternative to expensive physics-based simulations. We demonstrate the approach through a case study on power grid resilience planning against flood events in Texas, where we use 16 high-fidelity flood scenarios to generate 800 additional synthetic scenarios for validation. The results show that NORTA-generated scenarios accurately capture essential statistical properties, with the out-of-sample performance of first-stage decisions closely matching expectations from the original stochastic programming model. This framework enables decision-makers to assess the robustness of their solutions when obtaining additional real-world data is prohibitively expensive. The approach bridges machine learning and operations research by providing a practical solution to scenario generation challenges in stochastic programming.
\end{abstract}

\section{Introduction}
Stochastic programming has emerged as a powerful framework for decision-making under uncertainty across diverse domains including energy systems, supply chain management, finance, and transportation planning \cite{BirgLouv97}. The two-stage stochastic programming paradigm, where first-stage decisions are made before uncertainty realization followed by recourse actions, is particularly prevalent in practice. However, a critical challenge lies in scenario generation: obtaining representative samples of uncertain parameters often requires expensive simulations, field measurements, or expert assessments.

The intersection of machine learning and operations research offers promising solutions to this challenge. Generative models, which have greatly benefited domains like computer vision and natural language processing, can potentially address scenario scarcity in stochastic programming. These models learn the underlying distribution from limited data and generate synthetic samples that preserve essential statistical properties.

This work bridges ML and OR by demonstrating how the Normal-to-Anything (NORTA) generative model \cite{Cario97} can enhance the reliability of two-stage stochastic programming solutions through comprehensive out-of-sample testing. Unlike traditional approaches that rely on either computationally expensive physics-based simulations or oversimplified independent sampling, the NORTA model offers a middle ground: computationally efficient generation of scenarios that preserve both marginal distributions and correlation structures.

We make three key contributions: (1) we establish a framework for using generative models to validate stochastic programming solutions when scenario data is scarce; (2) we demonstrate that NORTA-generated scenarios accurately capture the statistical properties necessary for reliable out-of-sample testing; (3) we provide empirical evidence through a case study on power grid resilience planning, showing that decisions validated using NORTA-generated scenarios perform comparably to those validated with expensive physics-based simulations.

\section{Methods}
Consider a two-stage stochastic program:
\begin{subequations}
\begin{equation}
    \min_{x \in \mathcal{X}} \quad c^T x + \mathbb{E}_\xi[Q(x, \xi)],
\end{equation}
\begin{equation}
    \text{subject to:} \quad Ax = b, 
\end{equation}
\end{subequations}
where $x$ represents first-stage decisions, $\xi$ denotes the random vector, and $Q(x, \xi)$ is the recourse function. In practice, we typically solve the sample average approximation (SAA) version of this problem using $N$ scenarios $\{\xi^1, \ldots, \xi^N\}$:
\begin{subequations}
\begin{equation}
    \min_{x \in \mathcal{X}} \quad c^T x + \frac{1}{N}\sum_{i=1}^{N} Q(x, \xi^i)
\end{equation}
\begin{equation}
    \text{subject to:} \quad Ax = b.
\end{equation}
\end{subequations}
The real-world performance of the optimal first-stage solution determined by using the sample average approximation method depends critically on the scenarios used. When $N$ is small due to data scarcity, out-of-sample testing becomes essential to assess the robustness of the first-stage solution. However, generating additional scenarios for testing often faces the same cost constraints that limited the initial scenario set. To address this issue, we propose to generate synthetic scenarios using the limited available data. Specifically, using the $N$ available scenarios, we determine the empirical cumulative distribution function for each element in the random vector, thereby estimating its marginal distribution. Furthermore, we estimate the pairwise correlations between the elements of the random vector. These estimated parameters serve as inputs to the NORTA model. The NORTA model enables us to generate samples of multivariate random vectors with the specified marginal distributions and correlation structure. In the following sections, we provide details on how the parameters of the NORTA model are estimated and how we use them to generate additional scenarios for out-of-sample testing of the decisions obtained from the two-stage stochastic programming model.

\subsection{NORTA Modeling}

The NORTA model takes as input: (a) the marginal distributions $F_{X_i}$ for each variable $X_i$, and (b) the target $n \times n$ correlation matrix $\mathbf{\Sigma}_X$ of $\mathbf{X}$. As output, the NORTA model generates $n$-dimensional multivariate random vectors $\mathbf{X} = (X_1, X_2, \ldots, X_n)$ with marginal distributions and correlation structure as specified in the input.

\textbf{Parameter Estimation}: To generate samples with desired marginal distributions, the NORTA model transforms a base multivariate normal vector $\mathbf{Z} = (Z_1, Z_2, \ldots, Z_n)$ as follows:
\begin{equation}
\mathbf{X} = \begin{bmatrix}
    F_{X_1}^{-1}[\Phi(Z_1)] \\
    F_{X_2}^{-1}[\Phi(Z_2)] \\
    \vdots \\
    F_{X_n}^{-1}[\Phi(Z_n)]
\end{bmatrix}
\end{equation}
where $\Phi$ is the standard normal cumulative distribution function  and $F_{X_i}^{-1}$ is the inverse cumulative distribution function of $X_i$. To preserve the target correlations in $\mathbf{X}$, we must determine the correlation matrix $\boldsymbol{\Sigma}_Z$ for the base vector $\mathbf{Z}$. Interestingly, there exists a non-decreasing function $c_{ij}: \mathbb{R} \mapsto \mathbb{R}$ such that:
$$
    c_{ij}[\rho_Z(i,j)] =  \rho_X(i,j),
$$
where $\rho_Z(i,j)$ and $\rho_X(i,j)$ are elements of matrix $\boldsymbol{\Sigma}_Z$ and $\boldsymbol{\Sigma}_X$, respectively. The characteristics of the function $c_{ij}$ described in \cite{Cario97} enable us to estimate $\boldsymbol{\Sigma}_Z$ from $\boldsymbol{\Sigma}_X$ using a numerical line search method. The exact implementation of the same is presented in \cite{shukla_thesis}. 

\textbf{NORTA defectiveness}: The elements of matrix $\boldsymbol{\Sigma}_Z$ estimated using the line search method may not always yield a positive semi-definite matrix. As is discussed in \cite{ghosh_norta}, this problem becomes more severe as the dimension of random vector $\mathbf{X}$ increases. To address this issue, we solve a semi-definite program to find the closest positive semi-definite matrix to our estimated matrix and refer to it as $\boldsymbol{Y}$. This procedure ensures that the generated samples maintain the specified marginal distributions while approximating the target correlation structure, even in high-dimensional settings where exact correlation matching may be infeasible.

\textbf{Synthetic scenario generation}: Given the feasible correlation matrix $\boldsymbol{Y}$, we perform its Cholesky's decomposition $\boldsymbol{Y} = \boldsymbol{M}\boldsymbol{M}^T$ and generate samples using the algorithm described in Appendix \ref{sampling_algorithm}. 

\subsection{Out-of-sample testing}
We use the synthetic scenarios generated using the NORTA model to conduct out-of-sample testing. Specifically, we solve the two-stage program with limited number of true scenarios to get an optimal first stage decision $x^*$. Next, we fix this first stage decision $x^*$, and compute an estimate of the expected total cost of the two-stage stochastic program using $M$ scenarios as given by:
\begin{equation}
   \hat{v}_{OOS} = c^T x^* + \frac{1}{M}\sum_{j=1}^{M} Q(x^*, \tilde{\xi}^j),
\end{equation}
where $\tilde{\xi}^j$s are NORTA-generated scenarios. With $M$ scenarios where  ($M \gg N$), we can characterize the out-of-sample performance of the first-stage decisions. 

\section{Case study}

We demonstrate the effectiveness of our proposed approach through a case study examining transmission grid resilience decision-making against extreme flood events in the state of Texas. The transmission grid in the Texas Gulf Coast region experiences frequent flooding. To maintain operational continuity during flood events, protective structures are strategically positioned around critical power grid components to mitigate flooding impacts. To determine optimal locations for these protective structures, various stochastic programming models have been proposed \cite{Mohadese, 10287582, shukla_scenario-based_2021, SOUTO2022107545}. However, these models require flood scenarios as critical inputs, which are severely limited in availability and computationally expensive to generate. For our study, we obtained access to 16 high-fidelity flood scenarios referred to as MEOWs \cite{MEOW}. We use all available scenarios to determine optimal first-stage decisions regarding the placement of flood protection barriers and their required height thresholds to withstand flooding while ensuring uninterrupted grid operations. Leveraging these same scenarios, we generate 800 additional synthetic scenarios using the NORTA model. We evaluate the NORTA model's performance by assessing how closely the marginal distributions and pairwise correlations of the generated synthetic scenarios match the input marginal distributions and pairwise correlations estimated from the 16 high-fidelity flooding scenarios, which were generated through extensive flood simulations using a geoscience-based predictive flood model called SLOSH \cite{SLOSH}. For additional details on how we use the NORTA model for flood scenario generation, we refer to \cite{shukla_norta}. Finally, we compute the objective value of our two-stage model for each synthetic scenario and calculate various statistics to characterize the robustness of the first-stage decisions' performance. 

We provide additional details about the flood maps, transmission grid dataset, and our two-stage stochastic programming formulation for transmission grid resilience decision-making problem in Appendices \ref{meow_info}, \ref{grid_info}, and \ref{model_info}, respectively.

\section{Results and discussion}
We assess the NORTA model's accuracy in Section \ref{validation} and the out-of-sample performance of first-stage decisions from the two-stage stochastic programming model in Section \ref{oost}.

\subsection{Data validation of generated scenarios}\label{validation}
We evaluate the NORTA model's performance by examining how the 800 generated samples match the input marginal distributions and correlation matrix. For each marginal distribution, we calculate the Earth mover's distance (EMD) for univariate variable $X_i$:
\begin{equation}\label{emd}
    EMD_{i} = \int_{-\infty}^{\infty} |F_{X_i} - \hat{F}_{X_i}|,
\end{equation}
where $F_{X_i}$ is the input marginal and $\hat{F}_{X_i}$ is the estimated marginal from NORTA samples. We calculate EMD for each input marginal and report average, standard deviation, and quantiles in Table \ref{tab:norta_error}'s Earth mover's distance column. Similarly, the correlation error column shows statistics for differences between input pairwise correlations and estimated correlations from NORTA samples. The results demonstrate that the NORTA model maintain reasonable accuracy in generating correlated flood scenarios, even when the random vector $\mathbf{X}$ is 72-dimensional.

\subsection{Out-of-sample testing}\label{oost}
For the 800 scenarios that we generate using the NORTA model, we fix the first-stage decisions obtained by solving the two-stage model described in Appendix \ref{model_info} with 16 MEOW scenarios and solve the second-stage recourse problem to compute the load shed in each scenario. We repeat this process for each budget level from 0 to \$80M at \$10M intervals. Finally, we summarize the performance of the first-stage decisions on the out-of-sample data in Table \ref{tab:norta_oost}. Specifically, we present the load shed's average value, standard deviation, and quantiles across 800 scenarios for different budget levels. The table shows that the average performance of the decisions remains close to what we anticipate from the two-stage model's objective function.

\section{Conclusions}
This study proposes a framework for out-of-sample testing of two-stage stochastic programming models with limited sample availability. The framework uses available samples to construct empirical marginal distributions and correlation matrices of the random vector representing model uncertainty. This empirical information feeds into the NORTA model to generate additional samples that preserve key characteristics of the original limited data. We recommend using all limited data to determine optimal first-stage decisions while employing NORTA-generated samples to evaluate out-of-sample performance of these decisions. This approach helps decision-makers assess first-stage decision robustness against new but similar samples when obtaining additional data is costly. We demonstrate the framework's effectiveness through a case study on transmission grid resilience decision-making for flood scenarios in the Texas Gulf Coast region. 

While our framework demonstrates promise for enhancing out-of-sample testing in stochastic programming, several important limitations warrant careful consideration. First, the NORTA method exhibits degraded performance as the dimensionality of the random vector increases. Second, our validation approach focuses on preserving marginal distributions and pairwise correlations but does not guarantee accurate representation of the complete joint distribution, potentially missing higher-order dependencies, tail behaviors, and complex multivariate relationships that cannot be captured through pairwise correlations alone. This limitation introduces risk that out-of-sample performance estimates may systematically deviate from true performance. Therefore, we recommend that decision-makers interpret NORTA-based validation results as approximations rather than definitive performance guarantees, particularly in applications where tail risks or complex multivariate dependencies are critical to operational outcomes.

\bibliographystyle{plain}  
\bibliography{references}  

\appendix
\section{Sample generation using the NORTA model}\label{sampling_algorithm}
The following algorithm is used to generate synthetic samples for our case study.
\begin{algorithm}[h]
\caption{NORTA Sampling}
\begin{algorithmic}[1]
\STATE \textbf{Input:} Cholesky factor $M$, number of samples $m$
\STATE Initialize sample set $S = \emptyset$
\FOR{$i = 1$ to $m$}
    \STATE Generate $\hat{\mathbf{Z}} \sim \mathcal{N}(\mathbf{0}, \mathbf{I})$ (i.i.d. standard normal)
    \STATE Compute $\mathbf{Z} = M\hat{\mathbf{Z}}$
    \STATE Transform: $X_j = F_{X_j}^{-1}[\Phi(Z_j)]$ for $j = 1,\ldots,n$
    \STATE Add sample $\mathbf{X} = (X_1,\ldots,X_n)$ to $S$
\ENDFOR
\STATE \textbf{return} $S$
\end{algorithmic}
\end{algorithm}

\section{Flood maps}\label{meow_info}
We use NOAA's Maximum Envelopes of Water (MEOW) maps, derived from the SLOSH hurricane storm-surge model, to represent flooding scenarios in our two-stage stochastic program. SLOSH employs a parametric wind field model incorporating storm track, maximum wind radius, and pressure differential to simulate hurricane-induced storm surges. The MEOW composite product aggregates maximum surge values across multiple hurricane simulations with varying intensity, speed, direction, and tide levels, providing robust datasets for long-term mitigation planning while capturing forecast uncertainties. For computational tractability, we reduce the original 192 MEOW maps for the Texas coastal region to 16 scenarios by focusing on Category 5 hurricanes from four critical directions (west, west-north-west, north-west, and north-north-west) across four forward speeds, as these represent the most severe flooding conditions. Each scenario is treated as equally probable in our model, with flood heights at substations extracted from the corresponding SLOSH mesh cells and compared against hardening levels through constraints to determine the operational grid topology under different flooding conditions.

\section{Transmission grid}\label{grid_info}
To model the transmission grid, we use the ACTIVSg2000 dataset \cite{activsg2000}, which contains demand data for 2000 buses distributed across 1250 substations and connected via 3206 branches. The dataset was developed to preserve statistical similarities with Texas's actual power transmission grid while protecting the confidentiality of the actual grid details. We perform network reduction on the ACTIVSg2000 dataset to aggregate inland buses not exposed to flooding risks while preserving the detailed coastal grid topology that faces storm surge threats. This reduction maintains electrical equivalence with the original dataset while reducing the computation time of our two-stage stochastic programming model. To improve realism, we replace the synthetic coordinates with actual substation location coordinates from real-world infrastructure databases. This approach preserves the grid's electrical structure while accurately capturing flood risks for each substation based on their true geographic positions.

\section{Two-stage stochastic program}\label{model_info}
We model the transmission grid network for the Texas Gulf Coast region as a graph where nodes represent buses and edges represent branches connecting these buses. We assume that only substations (and therefore buses within them) are vulnerable to flooding impacts, as transmission lines are positioned well above ground level and are therefore considered immune to flooding. The two-stage stochastic program we consider operates as follows. In the first stage, decisions involve selecting which substations to protect with flood barriers and determining the flood height these barriers should withstand. After making these first-stage substation protection decisions, a flooding scenario is realized, causing unprotected substations or those with insufficient protection to fail. In the second stage, we solve a direct current-based power flow approximation model to determine optimal power routing through the remaining operational grid. We present a more detailed mathematical formulation of this two-stage stochastic programming model is the following subsections.

\subsection{Notation}
\noindent\textbf{Sets}
\allowdisplaybreaks
\begin{flalign*}
    \mathcal{I} &: \text{Set of substations indexed by $i$} && \\
    \mathcal{I}_f &: \text{Set of substations flooded in at least one scenario} && \\
    \mathcal{J} &: \text{Set of buses indexed by $j$} && \\
    \mathcal{B}_i &: \text{Set of buses at substation $i$} && \\
    \mathcal{K} &: \text{Set of scenarios indexed by $k$} && \\
    \mathcal{R} &: \text{Set of branches indexed by $r$} && \\
    \mathcal{N}_j^{in} &: \text{Set of branches incident on bus $j$ with power flowing into bus $j$} && \\
    \mathcal{N}_j^{out} &: \text{Set of branches incident on bus $j$ with power flowing out of bus $j$} &&
\end{flalign*}
\noindent\textbf{Parameters}
\begin{flalign*}
    M &: \text{An arbitrarily large constant} && \\
    I &: \text{Total investment budget for substation hardening} && \\
    f_i &: \text{Fixed cost of hardening at substation $i$} && \\
    v_i &: \text{Variable cost of hardening at substation $i$} && \\
    H_i &: \text{Maximum flood height to which substation $i$ can be hardened} && \\
    \Delta_{i}^k &: \text{Flood height at substation $i$ in scenario $k$ (a non-negative integer value)} && \\
    B_r &: \text{Susceptance of branch $r$} && \\
    F_r &: \text{Maximum power that can flow in branch $r$} && \\
    r.\mathrm{h}, r.\mathrm{t} &: \text{Head bus and tail bus of branch $r$} && \\
    D_j &: \text{Load at bus $j$} && \\
    \underline{G}_j, \overline{G}_j &: \text{Minimum and maximum generation at bus $j$} && \\
    \beta &: \text{Index of the reference bus} && \\
    p_k &: \text{Probability of scenario $k$} &&
\end{flalign*}

\noindent\textbf{Variables}
\allowdisplaybreaks
\begin{flalign*}
    y_i &: \text{Binary variable indicating whether substation $i$ is chosen for permanent hardening} && \\
    x_i &: \text{Non-negative integer variable indicating discrete height of hardening at substation $i$} && \\
    z_{j}^k &: \text{Binary variable indicating if bus $j$ is operational in scenario $k$} && \\
    s_{j}^k &: \text{Non-negative real variable, load satisfied at bus $j$ in scenario $k$} && \\
    g_{j}^k &: \text{Non-negative real variable, power generated at bus $j$ in scenario $k$} && \\
    u_{j}^k &: \text{Binary variable indicating if generator at bus $j$ is used in scenario $k$} && \\
    \alpha_{j}^k &: \text{Real variable indicating voltage phase angle of bus $j$ in scenario $k$} && \\
    e_{r}^k &: \text{Real variable indicating power flowing in branch $r$ in scenario $k$} &&
\end{flalign*}

\subsection{Mathematical formulation}
The two-stage stochastic programming model is given by:
\begin{equation}\label{SO_problem}
    \min_{x \in \mathcal{X}} \sum_{k \in \mathcal{K}} p_k\mathcal{L}(x,k),  
\end{equation}
where $\mathcal{X}$ represents the feasible set for the first-stage decisions and $\mathcal{L}(x,k)$ represents the load shed value when the first-stage decision is $x$ and flood scenario $k$ is realized. The feasible set $\mathcal{X}$ is defined by the constraints:
\begin{subequations}
\begin{equation}\label{budget_constraint}
        \sum_{i \in \mathcal{I}_f} f_i y_{i} + v_i x_{i} \leq I,
    \end{equation}
    \begin{equation}\label{constraint_tightening}
       x_{i} \leq H_i y_i, \quad \forall i \in \mathcal{I}_f. 
    \end{equation}
\end{subequations}
Specifically, constraint \eqref{budget_constraint} ensures that the total hardening cost does not exceed the available investment budget, and constraints \eqref{constraint_tightening} ensure that we can build a flood barrier no taller than $H_i$, and only if it is chosen for hardening.

The recourse function $\mathcal{L}(x,k)$ is defined as follows:
\begin{subequations}
\begin{equation}
    \mathcal{L}(x,k) = \min_{\textbf{$z^k$}, \textbf{$s^k$}, \textbf{$g^k$}, \textbf{$u^k$}, \textbf{$\alpha^k$}, \textbf{$e^k$}} \,\, \sum_{j \in \mathcal{J}} D_j - s_{j}^k,\label{objective}
\end{equation}
\allowdisplaybreaks
\begin{align}
    \text{subject to:} \quad & M(1-z_{j}^k) \geq \Delta_{i}^k - x_{i}, & \forall j \in \mathcal{B}_i, \quad \forall i \in \mathcal{I}_f, \label{linking_1} \\
    & 2Mz_{j}^k \geq 1 - 2(\Delta_{i}^k - x_{i}), & \forall j \in \mathcal{B}_i, \quad \forall i \in \mathcal{I}_f, \label{linking_2} \\
    & z_{j}^k = 1, & \forall j \in \mathcal{B}_i, \quad \forall i \in \mathcal{I} \setminus \mathcal{I}_f, \label{other_z} \\
    & s_{j}^k \leq D_jz_{j}^k, &  \quad  \forall j \in \mathcal{J},\label{supply_demand}\\
    & u_{j}^k \leq z_{j}^k, & \forall j \in \mathcal{J}, \label{flexible_generation}\\
    & u_{j}^k\underline{G}_j \leq g_{j}^k \leq u_{j}^k\overline{G}_j, & \forall j \in \mathcal{J},\label{capacity_constraint}\\
    & -z_{r.\mathrm{h}}^k F_r \leq e_{r}^k \leq z_{r.\mathrm{h}}^k F_r, & \forall r \in \mathcal{R},\label{edge_operational_1}\\
    & -z_{r.\mathrm{t}}^k F_r \leq e_{r}^k \leq z_{r.\mathrm{t}}^k F_r,  &  \forall r \in \mathcal{R},\label{edge_operational_2}\\
    & |e_{r}^k - B_r (\alpha_{r.\mathrm{h}}^k - \alpha_{r.\mathrm{t}}^k)|
                   \leq M(1 - z_{r.\mathrm{h}}^kz_{r.\mathrm{t}}^k), & \forall r \in \mathcal{R},\label{phase_angle_constraint}\\
    & \sum_{r \in \mathcal{N}_j^{out}}e_{r}^k - \sum_{r \in \mathcal{N}_j^{in}}e_{r}^k= g_{j}^k - s_{j}^k, & \forall j \in \mathcal{J}, \label{flow_balance}\\
    & -\pi \leq \alpha_{j}^k \leq \pi, & \forall j \in \mathcal{J}, \label{phase_side_constraint}\\
    & \alpha_{\beta }^k = 0. & \label{reference_bus}
\end{align}
\end{subequations}

In the aforementioned definition of the recourse function, the objective is to minimize the total load shed. Constraints \eqref{linking_1} and \eqref{linking_2} compare the flood height in scenario $k$ at each substation with the available protection level to determine whether the buses within it should be considered out of order. Constraints \eqref{other_z} ensure that all buses within substations that are not flooded remain operational. Constraints \eqref{supply_demand} ensure that the power consumed at bus $j$ cannot exceed the demand at that bus and that demand will be met only if that bus is operational. Constraints \eqref{flexible_generation} and \eqref{capacity_constraint} represent the generator dispatch decisions. Constraints \eqref{edge_operational_1} and \eqref{edge_operational_2} ensure that power can flow through a branch only if the buses at both ends are operational. Constraints \eqref{phase_angle_constraint} are phase angle constraints derived from the direct current-based power flow approximation. Constraints \eqref{flow_balance} are standard flow conservation constraints at each bus, constraints \eqref{phase_side_constraint} are side constraints on phase angles, and constraint \eqref{reference_bus} sets the phase angle of the reference or slack bus to zero.

\section{Tables}\label{table_section}
This section presents the complete numerical results supporting the validation and out-of-sample testing analyses discussed in Sections \ref{validation} and \ref{oost}. Table \ref{tab:norta_error} provides detailed statistics on the NORTA model's accuracy in preserving marginal distributions and correlation structures across the 72-dimensional flood scenario space. Table \ref{tab:norta_oost} presents comprehensive out-of-sample performance metrics for first-stage decisions across all budget levels, demonstrating the robustness of the stochastic programming solutions when validated against NORTA-generated synthetic scenarios.

\begin{table}[h]
    \centering
    \caption{The Earth mover's distance column shows the value of the mean, standard deviation, and distribution across various percentiles for the Earth mover's distance between the input marginal distribution and the distribution estimated from the NORTA samples across 72 flooded substations. The correlation error column shows the same metrics for the difference between pair-wise correlations of the input correlation matrix and the correlation matrix estimated from the NORTA samples.}
    \begin{tabular}{|l|c|c|}
    \hline
     &  Earth mover's distance& Correlation error\\
    \hline
    mean  &   0.073 &     0.041 \\
    std   &   0.038 &     0.038 \\
    min   &   0.001 &     0.000 \\
    25\%   &   0.052 &     0.013 \\
    50\%   &   0.071 &     0.030 \\
    75\%   &   0.103 &     0.058 \\
    max   &   0.146 &     0.329 \\
    
    \hline
    \end{tabular}
    \label{tab:norta_error}
\end{table}

\begin{table}[h]
    \centering
    \caption{The first row (SO estimate) represents the value of the expected load shed at different budget levels as was estimated from the two-stage stochastic programming (SO) model. The subsequent rows show the mean, standard deviation, and distribution of load shed across various percentiles in the second stage of the two-stage model at different budget levels (from 0 to \$80M). These values are computed using the 800 flood scenarios generated from the NORTA model for the out-of-sample testing of the first-stage substation hardening decisions for the Texas grid case study using the SO model.}
    \begin{tabular}{|l|c|c|c|c|c|c|c|c|c|}
    \hline
{} &      0  &      10 &      20 &      30 &      40 &      50 &      60 &      70 &    80 \\
\hline
 SO estimate & 4.462& 2.220& 1.374& 0.830& 0.479& 0.235& 0.060& 0.002&0.0\\
mean  &   4.475 &   2.231 &   1.381 &   0.832 &   0.481 &   0.237 &   0.060 &   0.002 &   0.0 \\
std   &   0.184 &   0.100 &   0.069 &   0.027 &   0.020 &   0.013 &   0.000 &   0.000 &   0.0 \\
min   &   3.980 &   1.952 &   1.216 &   0.774 &   0.447 &   0.215 &   0.059 &   0.002 &   0.0 \\
25\%   &   4.377 &   2.164 &   1.332 &   0.813 &   0.468 &   0.226 &   0.060 &   0.002 &   0.0 \\
50\%   &   4.481 &   2.225 &   1.388 &   0.831 &   0.478 &   0.238 &   0.060 &   0.002 &   0.0 \\
75\%   &   4.579 &   2.317 &   1.413 &   0.847 &   0.489 &   0.241 &   0.061 &   0.002 &   0.0 \\
max   &   4.975 &   2.445 &   1.523 &   0.894 &   0.531 &   0.264 &   0.061 &   0.002 &   0.0 \\
\hline
\end{tabular}
\label{tab:norta_oost}
\end{table}

\end{document}